\documentclass[a4paper,11pt]{amsart}
\usepackage{amssymb}
\usepackage{amsthm, amsmath, latexsym, amsfonts}
\usepackage[all]{xy}
\newtheorem{Lemma}{Lemma}
\newtheorem{Proposition}{Proposition}
\newtheorem{Theorem}{Theorem}
\newtheorem{remark}{{\em Remark}}
\newtheorem{Corollary}{Corollary}

\newtheorem{Question}{{\em Question}}

\newcommand{\Pic}{\mathrm{Pic}}

\newcommand{\Aut}{\mathrm{Aut}}

\newcommand{\Q}{\mathbb{Q}}

\renewcommand{\P}{\mathbb{P}}
\newcommand{\mg}{\mathcal{M}_g}
\newcommand{\mgn}{\mathcal{M}_{g,n}}

\newcommand{\mgbar}{\overline{\mathcal{M}}_g}
\newcommand{\mgnbar}{\overline{\mathcal{M}}_{g,n}}
\newcommand{\sgbar}{\overline{S}_g}
\newcommand{\sgnm}{\overline{S}_{g,n}}

\title{Moduli of Curves and Spin Structures \\ 
via Algebraic Geometry}
\author{Gilberto Bini and Claudio Fontanari}

\email{gilberto.bini@mat.unimi.it} \curraddr{Dipartimento di
Matematica, Universit\`a degli Studi di Milano \\ Via C. Saldini 50 \\
20133 Milano \\ Italy.}

\email{fontanar@science.unitn.it}\curraddr{
Universit\`a degli Studi di Trento \\
Dipartimento di Matematica \\
Via Sommarive 14 \\
38050 Povo (Trento) \\ Italy.}

\thanks{ {\em 2000 Mathematics Subject Classification}: 14H10, 14E08}
\begin{document}

\begin{abstract}
Here we investigate some birational properties of two collections of
moduli spaces, namely moduli spaces of
(pointed) stable curves and of (pointed) spin curves. In particular,
we focus on vanishings of Hodge numbers of type $(p,0)$ and on 
computations of the Kodaira dimension. 
Our methods are purely algebro-geometric and rely
on an induction argument on the number of marked points and the genus
of the curves (cf. \cite{ArbCor:98}).
\end{abstract}

\maketitle

\section{Introduction}

In the last decade, ideas from physics have reinvigorated the interest
in moduli spaces. As a consequence, the problem of investigating their
geometry has become more and more intriguing. Usually, such spaces
parametrize (pointed) genus $g$ complex curves with, possibly, extra
structures defined on $C$, such as a morphism to a fixed variety or a
spin structure. Moreover, they may be usually viewed as objects in
different categories, specifically as stacks, orbifolds, or
$\Q$-factorial projective schemes.

In the present paper, we focus on two classical moduli spaces,
namely $\mgnbar$ and $\sgnm$. We recall that the former one
parametrizes $n$-pointed genus $g$ stable curves and the latter one
parametrizes $n$-pointed quasi-stable curves $C$ of genus $g$ with a
line bundle whose second tensor power is isomorphic to the 
dualizing sheaf $\omega_C$: see Section \ref{spin} for precise 
definitions. Throughout, we regard them as projective normal varieties 
over the field of complex numbers.

It is impossible to condense in a short paragraph all the various
results and open problems in connection with these two varieties. We
refer the reader to, for instance, \cite{Corn}, \cite{Fab}, \cite{HM}
and the references cited therein. We just mention here that some
birational invariants of $\mgnbar$ are only conjecturally known, like
the Kodaira dimension of $\mgbar$, $17 \leq g \leq 23$, or partially
computed, like $(p,0)$-type Hodge numbers (cfr. \cite{Get1},
\cite{Get2}, \cite{GraPan:01}, and \cite{Pikaart:95}). Furthermore,
little is known about $\sgnm$. As a finite ramified covering of
$\mgnbar$, its geometry is even more complicated. Recent results on
the topology and the rational cohomology can be found, for example, in
\cite{Capo}, \cite{Cornalba:89}, \cite{Cornalba:91},
\cite{Harer:90}, \cite{Harer:93}. For a different, stack-theoretic approach
see \cite{AbrJar:01}, \cite{Jar1}, \cite{Jar2}.

As shown by Enrico Arbarello and Maurizio Cornalba in
\cite{ArbCor:98}, a purely algebro-geometric investigation of
$\mgnbar$ becomes more natural when the whole collection of
$\mgnbar$'s, and specific morphisms among them, is taken into account.
By an elementary double induction argument, it is thus possible to
prove various results on the rational cohomology of these spaces: see
\cite{Har1}, \cite{Har2} for related results which had been previously
known via other methods.

Here we apply a similar induction argument to $\mgnbar$ and
$\sgnm$. As a result, in Section~\ref{hodgenumbers} we derive some
vanishings of Hodge numbers of type $(p,0)$ for $\mgnbar$.
Next, in Section \ref{spin} we describe the rational Picard
group of $\sgnm$ and, as a corollary of \cite{Harer:93},
we determine generators and relations in the case $g \ge 9$, 
$n=0$. In Section~\ref{kodim}, we compute the Kodaira dimension of
spin moduli spaces in several cases and leave the other ones as open
questions, which we hope to address in the future with
different methods. To carry through some of these calculations, we
also compute the Kodaira dimension of ${\overline{\mathcal
M}}_{1,n}$. 
 
{\bf Acknowledgements.} The authors would like to thank M. Andreatta,
E. Ballico and S. Keel for useful remarks and suggestions.

\section{$(p,0)$ Hodge Numbers of Moduli Spaces of Stable Curves}
\label{hodgenumbers}

In \cite{ArbCor:98}, Enrico Arbarello and Maurizio Cornalba proved the
vanishing of some odd cohomology groups of $\mgnbar$. They applied an
inductive method that reduces the problem to checking such vanishings
for finitely many values of $g$ and $n$ in each odd degree
$k$. Unfortunately, if $k \ge 11$ the inductive machinery does not
work: indeed, it is well known (see, for instance, \cite{Pikaart:95},
proof of Corollary~4.7, or \cite{GraPan:01}, Proposition~2) that
$h^{11,0}(\overline{\mathcal{M}}_{1,11}) \ne 0$, where
$h^{p,q}(\mgnbar) = \dim H^{p,q}(\mgnbar)$ denote the Hodge numbers of
$\mgnbar$.  On the other hand, from the results of \cite{ArbCor:98} it
follows that $h^{p,0}(\mgnbar)=0$ for $p=1,3,5$. Here we complete the
picture by showing that there are no nonzero $p$ holomorphic
forms on $\mgnbar$ for $0 < p < 11$. Precisely, the following holds.

\begin{Theorem}\label{hodge}
Let $g$ and $n$ be non-negative integers, $n > 2-2g$.  If $0 < p < 11$,
then $h^{p,0}(\mgnbar)=0$.
\end{Theorem}

\proof As in \cite{ArbCor:98}, we use the long exact sequence of cohomology 
with compact supports:
$$
\ldots \to H^k_c(\mgn) \to H^k(\mgnbar) \to H^k(\partial \mgn) \to \ldots 
$$ 
and the vanishing of $H^k_c(\mgn)$ for
$$
k \le d(g,n) = \left \{  
\begin{array}{ll}
n-4  & \textrm{if} \hspace{0.1cm} g=0 \\
2g-2  & \textrm{if} \hspace{0.1cm} n=0 \\
2g-3+n  & \textrm{if} \hspace{0.1cm} g > 0, n > 0. 
\end{array}
\right.
$$
Moreover, since the morphism 
$$
H^k(\mgnbar) \to H^k(\partial \mgn)
$$
is compatible with the Hodge structures (see \cite{ArbCor:98}, p.~102), 
there is an injection 
$$
H^{p,0}(\mgnbar) \hookrightarrow H^{p,0}(\partial \mgn) 
$$
for $p \le d(g,n)$.

In order to conclude, it suffices to check that $h^{p,0}(\mgnbar)=0$ for 
every $p > d(g, n)$. This is in fact the base of the inductive procedure. 
When $g=0$, the space ${\overline{\mathcal M}}_{0,n}$ is rational \cite{Keel}; 
so there is nothing to check. If $g=1$, we are dealing with 
${\overline{\mathcal M}}_{1,n}$ for $n\leq 10$. 
As shown in \cite{Bel}, all these spaces are rational, so the base
of the induction holds true. Finally, when $g \geq 2$ and $p >
d(g,n)$, $\mgnbar$ is a unirational variety - see \cite{Logan:01},
Theorem~7.1). Hence, from a dominant rational map
$$
\P^{3g-3+n} \dashrightarrow \mgnbar
$$
we get an injective morphism
$$
H^{p,0}(\mgnbar) \hookrightarrow H^{p,0}(\P^{3g-3+n})
$$ as in \cite{GH}, p.~494, and we have $h^{p,0}(\P^{3g-3+n})=0$ for
every $p > 0$ (see for instance \cite{GH}, Corollary on p. 118). Thus
the claim is completely proved.

\qed

\section{The Rational Picard Group of Spin \\
Moduli Spaces}\label{spin}

The moduli space $S_g$ of smooth spin curves is a classical object,
which parametrizes pairs given by (smooth genus $g$ complex curve $C$,
and theta-characteristic on $C$). Since it is a non-trivial covering
of $\mg$, \emph{a priori} its geometry is much more complicated.
However, John Harer in \cite{Harer:90} and \cite{Harer:93} succeeded
in applying to $S_g$ his approach to the cohomology of moduli spaces
via geometric topology. Furthermore, in \cite{Cornalba:89} and
\cite{Cornalba:91}, Maurizio Cornalba constructed a geometrically
meaningful compactification $\sgbar$ of $S_g$. 
Here we determine $\Pic(\sgnm) \otimes \Q$ and give
explicit generators and relations. 
For reader's convenience we recall some basic definitions.

Let $X$ be a Deligne-Mumford semistable curve and let $E$ be a
complete, irreducible subcurve of $X$. $E$ is said to be
\emph{exceptional} when it is smooth, rational, and meets the other
components in exactly two points. Moreover, $X$ is said to be
\emph{quasi-stable} when any two distinct exceptional components of $C$
are disjoint. In the sequel, $\tilde{X}$ will denote the subcurve
$\overline{X \setminus \cup E_i}$ obtained from $X$ by removing all
exceptional components.

A \emph{spin curve} of genus $g$ (see \cite{Cornalba:89}, \S~2)
is the datum of a quasi-stable genus $g$ curve $X$ with an invertible 
sheaf $\zeta_X$ of degree $g-1$ on $X$ and a homomorphism of invertible 
sheaves
$$
\alpha_X: \zeta_X^{\otimes 2} \longrightarrow \omega_X
$$
such that 

(i) $\zeta_X$ has degree $1$ on every exceptional component of $X$;

(ii) $\alpha_X$ is not zero at a general point of every non-exceptional 
component of $X$.

\noindent Therefore, $\alpha_X$ vanishes identically 
on all exceptional components of $X$ and induces an isomorphism 
$$
\tilde{\alpha}_X: \zeta_X^{\otimes 2} \vert_{\tilde{X}} \longrightarrow 
\omega_{\tilde{X}}.
$$ 
In particular, when $X$ is smooth, $\zeta_X$ is just a theta-characteristic 
on $X$. 

\noindent Two spin curves 
$(X, \zeta_X, \alpha_X)$ and $(X', \zeta_{X'}, \alpha_{X'})$ are 
\emph{isomorphic} if there are isomorphisms $\sigma: X \to X'$ and 
$\tau: \sigma^*(\zeta_X') \to \zeta_X$
such that $\tau$ is compatible with the natural 
isomorphism between $\sigma^*(\omega_{X'})$ and $\omega_X$.

\noindent A \emph{family of spin curves} is a flat family of quasi-stable 
curves $f: \mathcal{X} \to S$ with an invertible sheaf $\zeta_f$ on 
$\mathcal{X}$ and a homomorphism
$$
\alpha_f: \zeta_f^{\otimes 2} \longrightarrow \omega_f
$$
such that the restriction of these data to any fiber of $f$ gives rise 
to a spin curve.

\noindent Two families of spin curves $f: \mathcal{X} \to S$ and 
$f': \mathcal{X}' \to S$ are \emph{isomorphic} if there are isomorphisms
$\sigma: \mathcal{X} \longrightarrow \mathcal{X}'$ and 
$\tau: \sigma^*(\zeta_{f'}) \longrightarrow \zeta_f$
such that $f = f' \circ \sigma$ and $\tau$ is compatible with the natural 
isomorphism between $\sigma^*(\omega_{f'})$ and $\omega_f$. 

\noindent Let $\sgbar$ be the set of isomorphism classes of spin curves of 
genus $g$ and $S_g$ be the subset consisting of classes of smooth curves. 
One can define a natural structure of analytic variety on $\sgbar$ 
(see \cite{Cornalba:89}, \S 5) in such a way that for any spin curve $X$ 
there is a neighbourhood of $[X]$ in $\sgbar$ and an isomorphism:
$$
U \cong B_X / \Aut(X),
$$ where $B_X$ is a $3g-3$-dimensional polydisk and $\Aut(X)$, the
group of automorphisms of the spin curve $X$, is a finite group (see
\cite{Cornalba:89}, Lemma~(2.2)). These spaces can be slightly generalized 
as follows:

\begin{eqnarray*}
\sgnm &:=& \{ [(C, p_1, \ldots, p_n; \zeta; \alpha)]: (C, p_1, \ldots,
p_n) \textrm{ is a genus $g$} \\ & & \textrm{quasi-stable projective
curve with $n$ marked points}; \\ & &\zeta \textrm{ is a line bundle
of degree $g-1$ on $C$ having} \\ & &
\textrm{degree $1$ on every exceptional component of $C$, and} \\ &
&\alpha: \zeta^{\otimes 2} \to \omega_C \textrm{
is a homomorphism which is} \\ & & \textrm{not zero at a general point
of every non-exceptional}\\ & & \textrm{component of $C$} \}.
\end{eqnarray*}

Analogously to $\sgbar$, these spaces are normal projective varieties
of complex dimension $3g-3+n$ with finite
quotient singularities. We point out the following fact:

\begin{Lemma}\label{iso}
Let $\Pic(\sgnm) := H^1(\sgnm, \mathcal{O}^*)$. There is a 
natural isomorphism 
$$ \Pic(\sgnm) \otimes \Q \stackrel{\cong}{\longrightarrow}
A_{3g-4+n}(\sgnm) \otimes \Q.
$$
\end{Lemma}

\proof Since $\sgnm$ is normal (see \cite{Cornalba:89}, Proposition~(5.2)),
there is an injection:
$$
\Pic(\sgnm) \hookrightarrow A_{3g-4+n}(\sgnm).
$$ Moreover, from the construction of $\sgnm$ it follows that the
singularities of $\sgnm$ are of finite quotient type, so every Weil
divisor is $\Q$-Cartier and there is a surjective morphism:
$$
\Pic(\sgnm) \otimes \Q \twoheadrightarrow A_{3g-4+n}(\sgnm) \otimes \Q.
$$
Hence the claim follows.
\qed

When $n=0$, it is also possible to give a more precise description of
of $\Pic(\sgbar) \otimes \Q$. 
To this end, recall that $\sgbar$ is the disjoint union of two irreducible
subvarieties $\sgbar^+$ and $\sgbar^-$ which consist of the even and
the odd spin curves of genus $g$ (see \cite{Cornalba:89},
Lemma~(6.3)), respectively. The following crucial result was obtained
by John Harer via geometric topology (see \cite{Harer:93},
Corollary~1.3):

\begin{Theorem}\label{harer}
Let $\mg[\varepsilon]$ denote either $\sgbar^+ \cap S_g$ or 
$\sgbar^- \cap S_g$. Then $\Pic(\mg[\varepsilon]) \otimes \Q := 
H^1(\mg[\varepsilon], \mathcal{O}^*) \otimes \Q$ has rank $1$ for 
$g \ge 9$. 
\end{Theorem}
 
For any family $f: \mathcal{X} \to S$ of spin curves, 
$M_f := \det Rf_* \zeta_f$ is a line bundle on $S$. Let $M$ denote the 
corresponding line bundle on $\sgbar$ associated to the universal 
family on $\sgbar$ (as usual, notice that for $g \ge 4$ the locus 
of spin curves with automorphisms has complex codimension $\ge 2$). 
Let $\mu^+$ (resp. $\mu^-$) be the class of $M$ in $\Pic(\sgbar^+)$ 
(resp. $\Pic(\sgbar^-)$.
The boundary $\partial \sgbar = \sgbar \setminus S_g$ is the union of 
irreducible components $A_i^+$, $B_i^+$ (contained in $\sgbar^+$) and 
$A_i^-$, $B_i^-$ (contained in $\sgbar^-$), which are completely 
described in \cite{Cornalba:89}, \S~7. Let $\alpha_i^+$, $\beta_i^+$, 
$\alpha_i^-$, $\beta_i^-$ denote the corresponding classes in 
$A_{3g-4}(\sgbar)$. The following result is contained in 
\cite{Cornalba:89}, Proposition~(7.2):
\begin{Proposition}\label{independent}
If $g > 2$ is odd, the classes $\mu^+$, $\mu^-$, $\alpha_i^+$, 
$\alpha_i^-$, $\beta_i^+$, $\beta_i^-$, $i= 0, \ldots, (g-1) / 2$, 
are independent. 
If $g > 2$ is even, the same is true of the classes $\mu^+$, $\mu^-$, 
$\alpha_i^+$, $\alpha_i^-$, $\beta_i^+$, $i=0, \ldots, g/2$, and 
$\beta_i^-$, $i= 0, \ldots, g / 2 - 1$.
\end{Proposition}

Now we are ready to state and prove a description by generators 
and relations of the rational Picard group of $\sgbar$:

\begin{Corollary}\label{pic}
Assume $g \ge 9$. If $g$ is odd, then $\Pic(\sgbar^+)$ (resp. 
$\Pic(\sgbar^-)$) is freely generated over $\Q$ by the classes 
$\mu^+$, $\alpha_i^+$, $\beta_i^+$, $i= 0, \ldots, (g-1) / 2$
(resp. by the classes $\mu^-$, $\alpha_i^-$, $\beta_i^-$, 
$i= 0, \ldots, (g-1) / 2$). 
If $g $ is even,  then $\Pic(\sgbar^+)$ (resp. $\Pic(\sgbar^-)$) is freely 
generated over $\Q$ by the classes $\mu^+$, $\alpha_i^+$, $\beta_i^+$, 
$i= 0, \ldots, g / 2$ (resp. by the classes $\mu^-$, $\alpha_i^-$, 
$i= 0, \ldots, g / 2$, and $\beta_i^-$, $i= 0, \ldots, g / 2 - 1$).
\end{Corollary}

\proof By Lemma~\ref{iso}, we may use the exact sequence
$$
A_{3g-4}(\sgbar \setminus S_g) \to A_{3g-4}(\sgbar) \to A_{3g-4}(S_g)
\to 0
$$
to conclude that $\Pic(\sgbar) \otimes \Q$ is generated by 
the generators of $A_{3g-4}(S_g)$ together with the set of boundary classes 
of $\sgbar$. From Theorem~\ref{harer} it follows that $\Pic(S_g) \otimes \Q$ 
is generated by the classes $\mu^+$ and $\mu^-$, therefore $\Pic(\sgbar) 
\otimes \Q$ is generated by the classes $\mu^+$, $\mu^-$, $\alpha_i^+$, 
$\alpha_i^-$, $\beta_i^+$, $\beta_i^-$.
By Proposition~\ref{independent}, all these classes are independent, so 
Corollary~\ref{pic} is proved. 

\qed

\section{Calculation of Kodaira Dimensions}\label{kodim}

In this section, we calculate the Kodaira dimension of some spin 
moduli spaces. We recall that the Kodaira dimension is an important
birational invariant in the classification of projective varieties. As
a first general result, we prove the following

\begin{Proposition} Fix $g$ and $n$ non-negative integers, $n >2-2g$. 
Any irreducible component of $\sgnm$ is of general type whenever
$\mgnbar$ is.
\end{Proposition}

\proof Since any component of $\sgnm$ is a non-trivial ramified
covering of $\mgnbar$, then the claim follows from \cite{Ueno:75},
Theorem~6.10. \qed

In particular, $\sgbar$ is of general type for $g \geq 24$. When $n
\geq 1$, Logan determines in \cite{Logan:01} integers
$\widetilde{n}(g)$, $g \geq 2$, such that $\mgnbar$ is of general type
when $n \geq \widetilde{n}(g)$. As pointed out in Section \ref{spin},
when $g=0$, spin moduli spaces are rational since they are isomorphic
to ${\overline {\mathcal M}}_{0,n}$. To tackle the genus one case, we
first need to compute the Kodaira dimension of ${\overline{\mathcal
M}}_{1,n}$, which does not seem to be thoroughly dealt with in the
literature. It turns out that $\kappa({\overline{\mathcal M}}_{1,n})$
varies with the number of marked points.

As proved in \cite{Bel}, the moduli space of $n$-pointed genus $1$
curves is rational for $n \leq 10$, hence $\kappa({\overline{\mathcal
M}}_{1,n})= - \infty$. To compute the Kodaira dimension for $n \geq
11$, we first need to express the canonical divisor
$K_{{\overline{\mathcal M}}_{1,n}}$ in terms of generators of the
rational Picard group of ${\overline{\mathcal M}}_{1,n}$. We briefly
recall such generators and some of their relations: for more details
the reader is referred, for instance, to \cite{ArbCor:87}.

As usual, we denote by $\lambda$ the first Chern class of the Hodge
bundle whose fiber over the element $[C;p_1, \ldots, p_n]\in
{\overline{\mathcal M}}_{1,n}$ is $H^0(C, \omega_C)$, where $\omega_C$
is the dualizing sheaf of $C$. Next, we denote by $\psi_i$, $1\leq i
\leq n$, the first Chern class of the line bundle whose fiber over
$[C;p_1, \ldots, p_n]$ is the cotangent space of $C$ at the (smooth)
point $p_i$. Finally, we denote by $\delta_{irr}$ and $\delta_{0,S}$
the classes corresponding to boundary divisors. Here $\delta_{irr}$ is
the (rational) Poincar\'{e} dual of the locus of curves with one
non-disconnecting node and $n$ marked points. The class $\delta_{0,S}$, 
$|S|\geq 2$, corresponds the locus of curves with a disconnecting node
whose removal creates two connected components: one of genus $0$ with
the marked points labelled by the elements of $S$ and the other one of
genus $1$ with the marked points labelled by the elements of
$S^c$. Last, note that by $\kappa_m$ and $\widetilde{\kappa}_m$,
$m\geq 1$ we mean the Mumford classes as described in
\cite{ArbCor:98}. 

We finally recall that $\Pic(\mgnbar)\otimes \Q$ is isomorphic to the
Picard group of the moduli stack of $n$-pointed genus $g$ stable
curves. We shall first compute the canonical class of this stack so to
deduce $K_{{\overline{\mathcal M}}_{1,n}}$.

\begin{Proposition}
\label{canonico1}
For any non-negative integer $ n \geq 3$, 
\begin{eqnarray}
\label{expcan}
K_{{\overline{\mathcal M}}_{1,n}} &= &(n-11)\lambda +
(n-3)\delta_{0, \{1, \ldots, n\}} \\ &+& \sum_{\underset{S
\subset\{1,\ldots, n\}}{|S| \geq 2}, S\neq \{1,\ldots,
n\}}\left(|S|-2\right)\delta_{0,S}.
\end{eqnarray}
\end{Proposition}
\proof Suppose $\rho: {\mathcal C} \rightarrow B$ is a family of
$n$-pointed genus $1$ curves over a smooth base $B$ and with general
smooth fiber. Let $\sigma_1, \ldots, \sigma_n$ denote the $n$
canonical sections and define $D$ to be the sum of the divisors
corresponding to them. As in \cite{HM}, we apply
Grothendieck-Riemann-Roch Theorem to the bundle $\Omega^1_{{\mathcal
C}/B}(D) \otimes \omega_{{\mathcal C}/B}$, where $\Omega^1_{{\mathcal
C}/B}$ is the sheaf of relative K\"{a}hler differentials and
$\omega_{{\mathcal C}/B}$ is the relative dualizing sheaf of
$\rho$. Therefore,
\begin{eqnarray}
\label{canonicstack}
K&=& c_1(\rho_{*}(\Omega^1_{{\mathcal C}/B}(D) \otimes
\omega_{{\mathcal C}/B})) \\ \nonumber &=&
\rho_{*}\left(Td^{\vee}_2(\Omega^1_{{\mathcal C}/B}) +
ch_2(\Omega^1_{{\mathcal C}/B}(D) \otimes \omega_{{\mathcal C}/B})\right) \\
&+& \rho_{*}\left(Td^{\vee}_1(\Omega^1_{{\mathcal C}/B})ch_1(\Omega^1_
{{\mathcal C}/B}(D) \otimes \omega_{{\mathcal C}/B})\right), \nonumber
\end{eqnarray}
where $Td_j^{\vee}$ and $ch_j$ denote the degree $j$ term of the Todd
class $Td^{\vee}$ and of the Chern character $ch$, respectively.

If $\eta$ is the class of the locus of nodes of fibers of
${\mathcal C}$ over $B$, we have
$$
Td^{\vee}_1(\Omega^1_{{\mathcal C}/B})= -\frac12c_1(\omega_{{\mathcal
C}/B}),
$$
$$
Td^{\vee}_2(\Omega^1_{{\mathcal C}/B})=
\frac{1}{12}\left(c_1^2(\omega_{{\mathcal C}/B}) + \eta\right).
$$

Analogously,
$$
ch_1(\Omega^1_{{\mathcal C}/B}(D)\otimes \omega_{{\mathcal C}/B})=
c_1(\omega_{{\mathcal C}/B}(D))+ c_1(\omega_{{\mathcal C}/B}),
$$
\begin{eqnarray*}
ch_2(\Omega^1_{{\mathcal C}/B}(D) \otimes \omega_{{\mathcal C}/B}))&=&
\frac{1}{2}c_1^2(\omega_{{\mathcal C}/B}) -\eta \\ &+&
c_1(\omega_{{\mathcal C}/B}(D))c_1(\omega_{{\mathcal C}/B}) + \frac12
c_1^2(\omega_{{\mathcal C}/B}(D)).
\end{eqnarray*}

By the definition of $\kappa_1$ and ${\widetilde{\kappa}}_1$, we
have
\begin{eqnarray}
\label{ameta}
K&=&\frac{1}{12}{\widetilde{\kappa}}_1+\frac{1}{2}\rho_{*}\left(c_1(\omega_
{{\mathcal C}/B}(D))c_1(\omega_{{\mathcal C}/B})\right) \nonumber \\
&-& \frac{11}{12}\rho_{*}(\eta) + \frac{1}{2}\kappa_1.
\end{eqnarray} 

On the other hand, by \cite{ArbCor:96}, 
$$
\kappa_1 = {\widetilde{\kappa}}_1 + \sum_{i=1}^n \psi_i,
$$
and
$$
\frac{1}{2}\rho_{*}\left(c_1(\omega_{{\mathcal
C}/B}(D))c_1(\omega_{{\mathcal C}/B})\right)= \frac{1}{2} \widetilde{\kappa}_1 
+ \sum_{i=1}^n\psi_i.
$$

Therefore, by Example 2.1 in \cite{Bin:03}, we have
\begin{equation}
\label{canigen}
K=13\lambda +\sum_{i=1}^n\psi_i -2\delta,
\end{equation}
where 
$$ 
\delta:=\delta_{irr} +
\sum_{\underset{S\subset\{1,\ldots,n\}}{|S|\geq 2}}\delta_{0,S}.
$$ 

Moreover, in genus $1$ (see \cite{ArbCor:98}), this can be rewritten
as
$$
K=(n-11)\lambda +\sum_{\underset{S \subset\{1,\ldots, n\}}{|S| \geq
2}}\left(|S|-2\right)\delta_{0,S}.
$$

Since the map from the moduli stack of $n$-pointed, $n \geq 3$, genus
$1$ curves to the (coarse) moduli space is ramified along the divisor
$\delta_{0, \{1, \ldots, n\}}$, the claim follows.

\qed

\begin{remark} {\rm When $g=1$ and $n=1$, the (coarse) moduli space is
isomorphic to ${\mathbb P}^1$, so the canonical class is known. When
$g=1$, and $n=2$, then (\ref{canigen}) simplifies to
$-9\lambda$. Analogously to the case $n \geq 3$, the canonical class
is $-9\lambda - \delta_{0, \{1,2\}}$.}
\end{remark}

\begin{remark}
{\rm In \cite{Logan:01}, a formula for the canonical divisor of
$\mgnbar$ is given. The proof relies on the corresponding formula for
$K_{{\overline{\mathcal M}}_g}$, the canonical divisor of
${\overline{\mathcal M}}_g$, obtained in \cite{HarM} only for $g \geq
2$ via Grothendieck-Riemann-Roch Theorem. Up to the authors'
knowledge, an analogous formula in genus $1$ is not explicitely stated
in the literature.}
\end{remark}

We can now complete the computation of $\kappa({\overline{\mathcal
M}}_{1,n})$ for each $n \geq 1$.  In fact, the following holds.

\begin{Theorem}\label{kodaira1} We have
$$
\kappa({\overline{\mathcal M}}_{1,n})=
\left\{
\begin{array}{ccc}
0 & & n = 11, \\
1 & & n \ge 12.
\end{array}
\right.
$$
\end{Theorem}

\proof 
By Proposition \ref{canonico1}, $K_{{\overline{\mathcal M}}_{1,11}}$ 
is an effective divisor, hence $\kappa({\overline{\mathcal M}}_{1,11})\ge 0$.
On the other hand, ${\mathcal M}_{1,11}$ is birational to a hypersurface 
$X$ of $({\mathbb P}^2)^6$ of multidegree $(3, \ldots, 3)$
(see Remark 1.2.4 in \cite{Bel}). By adjunction, we obtain that 
$K_X$ is trivial; in order to compute the Kodaira dimension of $X$, 
let $f: Y \to X$ be the normalization map. We have $K_Y = f^* K_X - 
\Delta$, where the conductor $\Delta$ is an effective divisor. 
It follows that $\kappa({\overline{\mathcal M}}_{1,11}) = 
\kappa(X) = \kappa(Y) \le 0$, so the case $n=11$ is over.   

Next, again by Proposition \ref{canonico1}, $K_{{\overline{\mathcal
M}}_{1,n}}$, $ n \geq 12$, is the sum of two effective divisors, i.e.,
$L:=(n-11)\lambda$ and 
$$
E:=(n-3)\delta_{0, \{1, \ldots, n\}}+\sum_{\underset{S\subset\{1,\ldots,
n\}}{|S| \geq 2}, S\neq \{1,\ldots, n\}}\left(|S|-2\right)\delta_{0,S}.
$$
Therefore, the Kodaira dimension of ${\overline{\mathcal M}}_{1,n}$ is
greater than or equal to the Iitaka dimension of the divisor $L$. On
the other hand, let
\begin{equation}\label{forgetful}
\pi: {\overline{\mathcal M}}_{1,n} \rightarrow {\overline{\mathcal M}}_{1,1}
\end{equation}
be the morphism which forgets the last $n-1$ points and passes to the
stable model. 
If $\lambda_1$ denotes the first Chern class of the Hodge bundle on
${\overline{\mathcal M}}_{1,1}$, then we have $\lambda =
\pi^*(\lambda_1)$ (see for instance \cite{ArbCor:87}, (6)). Moreover, 
since $\lambda_1$ is ample on ${\overline{\mathcal M}}_{1,1}$, we obtain
$$
\kappa({\overline{\mathcal M}}_{1,n}, L) \ge \kappa({\overline{\mathcal
M}}_{1,1}, (n-11)\lambda_1)=1.
$$
This proves that $\kappa({\overline{\mathcal M}}_{1,n}) \geq
1$. However, the fiber of $[C;p] \in {\overline{\mathcal M}}_{1,n}$
under $\pi$ can be viewed as the quotient by a finite group of an open
Zariski subset of the product 
$$
\underbrace{C \times \ldots \times C}_{\hbox{ (n-1) times}},
$$
hence $\kappa(\pi^{-1}([C;p])) \le 0$.
By Theorem~6.12 in \cite{Ueno:75}, it follows that 
$$
\kappa({\overline{\mathcal M}}_{1,n}) \leq 
\kappa(\pi^{-1}([C;p])) + \dim ({\overline{\mathcal M}}_{1,1}) \le 1.
$$
Thus the claim is completely proved.  

\qed

\begin{Corollary}
For any $n \geq 1$, ${\overline{\mathcal M}}_{1,n}$ is never of general type.
\end{Corollary}

\bigskip

Next, we turn to moduli spaces of pointed spin curves of genus $1$. 
Recall from Section~\ref{spin} that 
$\overline{S}_{1,n}$
is the compactification \emph{\`a la} Deligne-Mumford of the 
moduli space of $n$-pointed smooth elliptic curves with a 
theta-characteristic. Notice that 
$\overline{S}_{1,n}$ is the disjoint 
union of $\overline{S}_{1,n}^+$
and $\overline{S}_{1,n}^-$, which 
correspond to even and odd theta-characteristics, respectively.
However, since over an elliptic curve there is only one odd 
theta-characteristic (namely, the structural sheaf), there is 
a natural isomorphism $\overline{S}_{1,n}^-
\cong \overline{\mathcal M}_{1,n}$. 
>From now onwards, we thus focus our attention on 
$\overline{S}_{1,n}^+$. 

In order to prove that 
$\kappa(\overline{S}_{1,n}^+) = - \infty$ 
for $n \le 10$, we are going to show that 
$\overline{S}_{1,n}^+$ is uniruled 
whenever $\overline{\mathcal M}_{1,n}$ is. 
Indeed, the following holds:

\begin{Lemma}
Let $p: \overline{S}_{1,n}^+ \to
\overline{\mathcal M}_{1,n}$ denote the natural projection.
If $C$ is a rational curve in $\overline{\mathcal M}_{1,n}$, 
then there exists a rational curve $D$ in
$\overline{S}_{1,n}^+$ such that 
$p(D)=C$. 
\end{Lemma} 

\proof The proof is by induction on $n$. 

If $n=1$, then we have $\overline{\mathcal M}_{1,n} \cong \P^1$ and 
$\overline{S}_{1,n}^+ \cong \P^1$,
so in this case the stated property is obvious.

Assume, now, $n > 1$. It is easy to check that 
$$
\overline{S}_{1,n}^+ =
\overline{\mathcal M}_{1,n} \times_{\overline{\mathcal M}_{1,n-1}}
\overline{S}_{1,n-1}^+.
$$
Therefore, we have the following commutative diagram:

\begin{equation*}
\xymatrix{
\P^1 \ar@/_/[ddr]_g \ar[dr] \ar@/^/[drr]^f \\
& \overline{S}_{1,n}^+ \ar[d] \ar[r]
               & \overline{\mathcal M}_{1,n} \ar[d]^\pi             \\
& \overline{S}_{1,n-1}^+ \ar[r] 
& \overline{\mathcal M}_{1,n-1}}
\end{equation*}
\vskip .1cm
\centerline{
Diagram 1: Uniruledness of spin moduli spaces for $g=1, n\leq 10$.}

\vskip .2cm
\noindent In Diagram 1, $f$ exists by hypothesis and $g$ exists by the
inductive assumption (just notice that $\pi \circ f (\P^1)$ is not a
point since the fibers of $\pi$ does not contain rational curves).
Hence the claim follows from the universal property of the fibered
product.

\qed

Finally, we consider the case $n \ge 12$.

\begin{Proposition}
\label{kod1s}
Let $n \ge 12$ be a non-negative integer. Then the Kodaira dimension of 
$\overline{S}_{1,n}^+$ is $1$.
\end{Proposition}

\proof Since the natural projection 
$$
p: \overline{S}_{1,n}^+ \longrightarrow
\overline{\mathcal M}_{1,n}
$$
is a surjective map between normal varieties of the same dimension,
it follows that 
$$
\kappa(\overline{S}_{1,n}^+) \ge 
\kappa(\overline{\mathcal M}_{1,n}) = 1
$$
(see \cite{Ueno:75}, Theorem~6.10). 

On the other hand, the fiber of the forgetful map
$$
\overline{S}_{1,n}^+ \longrightarrow
\overline{S}_{1,1}^+
$$
is precisely the same as that of the morphism (\ref{forgetful}). 
Hence, exactly as in the proof of Theorem~\ref{kodaira1}, 
we can deduce that
$$
\kappa(\overline{S}_{1,n}^+) \le 
0 + \dim(\overline{S}_{1,1}^+) \le 1.
$$
Thus the proof follows.

\qed

By the same arguments as in Proposition \ref{kod1s}, we get the
estimate
$$
0 \leq \kappa(\overline{S}_{1,11}^+) \leq 1.
$$

Another reason why $\kappa(\overline{S}_{1,11}^+) \geq 0$ is given by 
the following result.

\begin{Proposition}
There exists a nonzero holomorphic form of degree $11$ on
$\overline{S}_{1,n}^+$.
\end{Proposition}

\proof Since the natural projection 
$$
p: \overline{S}_{1,n}^+ \longrightarrow
\overline{\mathcal M}_{1,n}
$$
is a surjective map between normal varieties of the same dimension,
the induced map
$$ p^*:H^{11,0}(\overline{\mathcal M}_{1,n}, \Q) \rightarrow
H^{11,0}(\overline{S}_{1,n}^+, \Q)
$$
is injective.

\qed

\bigskip

We end this section with a couple of natural questions. 
 
\begin{Question} {\rm Is the Kodaira dimension of
$\overline{S}_{1,11}^+$ zero\,?}
\end{Question}

\begin{Question} {\rm Is any irreducible component of $\sgnm$ unirational
whenever the corresponding moduli space $\mgnbar$ is\,?}
\end{Question}

\end{document}